\documentclass[12pt]{amsart}
\usepackage{amscd,times,fullpage}
\usepackage[utf8x]{inputenc}

\usepackage{amsthm,amsmath}
\usepackage{amssymb}
\usepackage{graphicx}
\usepackage{amsthm}
\usepackage{enumerate}
\usepackage{dsfont}
\usepackage{mathtools}
\usepackage{tikz-cd}
\usetikzlibrary{cd}
\usetikzlibrary{positioning}
\usetikzlibrary{matrix}
\usepackage[hidelinks]{hyperref}

\hypersetup {
	pdfpagemode = {UseNone},
	pdftitle = {On weak and strict relatives K\"ahler manifolds},
	pdfauthor = {Giovanni Placini},
	pdflang = {en-UK}
} 

\newcommand{\R}{\mathbb{R}}

\newcommand{\Z}{\mathbb{Z}}

\newcommand{\CP}{\mathbb{C}\mathrm{P}}

\newcommand{\CH}{\mathbb{C}\mathrm{H}}

\newcommand{\Hol}{\mathrm{Hol}}
\newcommand{\hol}{\mathfrak{hol}}

\newcommand{\C}{\mathbb{C}}            
\newcommand{\de}{\partial}

\newcommand{\K}{K\"{a}hler }

\newcommand{\ov}[1]{\overline{#1}}

\newcommand{\deb}{\ov\partial}

\newcommand{\Id}{\operatorname{Id}}

\newcommand{\so}{\mathfrak{so}}

\newcommand{\lra}{\longrightarrow}

\newtheorem{theor}{Theorem}

\newtheorem{defin}[theor]{Definition}
\newtheorem{lem}[theor]{Lemma}
\newtheorem{cor}[theor]{Corollary}

\newtheorem{ex}{Example}

\makeatletter
\newtheorem*{rep@theorem}{\rep@title}
\newcommand{\newreptheor}[2]{%
	\newenvironment{rep#1}[1]{%
		\def\rep@title{#2 \ref{##1}}%
		\begin{rep@theorem}}%
		{\end{rep@theorem}}}
\makeatother

\newreptheor{theor}{Theorem}

\usepackage{nicefrac}

\begin{document}

	\title[On weak and strict relatives K\"ahler manifolds]{On weak and strict relatives K\"ahler manifolds}

	\author{G.~Placini}
	\address{Dipartimento di Matematica e Informatica, Universit\'a degli studi di Cagliari, Via Ospedale 72, 09124 Cagliari, Italy}
	\email{giovanni.placini@unica.it}

	\date{\today ; {\copyright  G.~Placini 2023}}
	
	\subjclass[2010]{Primary 53C42; Secondary 53C24, 32Q15, 53B35.}
	\keywords{Relatives K\"ahler manifolds; holomorphic isometries, Weak relatives K\"ahler manifolds, K\"ahler submanifolds}
	\thanks{The author is supported by INdAM and  GNSAGA - Gruppo Nazionale per le Strutture Algebriche, Geometriche e le loro Applicazioni, by GOACT - Funded by Fondazione di Sardegna and funded the National Recovery and Resilience Plan (NRRP), Mission 4 Component 2 Investment 1.5 - Call for tender No.3277 published on December 30, 2021 by the Italian Ministry of University and Research (MUR) funded by the European Union – NextGenerationEU. Project Code ECS0000038 – Project Title eINS Ecosystem of Innovation for Next Generation Sardinia – CUP F53C22000430001- Grant Assignment Decree No. 1056 adopted on June 23, 2022 by (MUR)}

	\begin{abstract}
	We study \K manifolds that are (weak) relatives, that is, \K manifolds which share a (locally isometric) submanifold. In particular, we prove that if two \K manifolds are weak relatives and one of them is projective, then they are relatives. Moreover, we introduce the notion of strict relatives K\"ahler manifolds and provide several nontrivial examples.
	\end{abstract}
	
	\maketitle
	
	\section{Introduction and statements of the main results}\label{sectionint}
The study of holomorphic isometries of \K manifolds provides one of the most striking examples of rigidity in \K geometry.
For instance, in his seminal work \cite{Calabi53Isometric} Calabi proved that it does not exist a holomorphic isometry between finite dimensional complex space forms of different type.
After Calabi, the topic has received great interest and many weakening of the holomorphic isometry condition have been studied.
One of these generalizations is the following: which \K manifolds share a \K submanifold?
 This question was first considered and answered in the negative in the case of complex space forms by Umehara in \cite{Umehara87Relatives}. 
\K manifolds which share a \K submanifold were christened relatives by Di Scala and Loi \cite{DiscalaLoi10Relatives} who proposed their systematic study. In particular, they gave the following definition.
	\begin{defin}[\cite{DiscalaLoi10Relatives}]\label{DefRelatives}
	Two \K manifolds $M_1$ and $M_2$ are \textbf{relatives} (or are related) if there exists a \K manifold $X$ and two holomorphic isometries $\varphi_i:X\lra M_i$ for $i=1,2$.
\end{defin}
Observe further that, since we are not assuming the compactness of $X$, we can consider $\varphi_i$ to be a holomorphic isometric embedding by restricting to a open neighbourhood of a point in $X$.
Umehara \cite{Umehara87Relatives} proved that complex space forms with holomorphic sectional curvature of different sign are not relatives extending the aforementioned Calabi's result.
Moreover, in \cite{DiscalaLoi10Relatives} Di Scala and Loi proved the same result for a bounded domain and a projective \K manifold. 
The reader can refer to \cite{Yuan19SurveyRelatives} and references therein for a recent survey on this topic.
Even more recently, several generalizations of these results have been proven by Loi and Mossa \cite{LoiMossa22HBD,LoiMossa22BlowUp,LoiMossa23Rigidity},  and Cheng and Hao \cite{ChengHao21Relatives}.
This rigidity phenomenon is peculiar of \K manifolds as very often the aforementioned results rely heavily on the compatibility of the Riemannian metric and the complex structure. 
It is then natural to weaken this notion and study which pairs of \K manifolds share a Riemannian submanifold which is not necessarily a holomorphic submanifold.
	\begin{defin}[\cite{DiscalaLoi10Relatives}]\label{DefWeaklyRelatives}
	Two \K manifolds $M_1$ and $M_2$ are \textbf{weak relatives} (or are weakly related) if there exists two locally isometric \K manifolds $X_1$ and $X_2$ of complex dimension$\geq 2$ and holomorphic isometries $\varphi_i:X_i\lra M_i$ for $i=1,2$.
\end{defin}
Observe that the condition on the dimension of the submanifolds is necessary. In fact, any isometry between two one dimensional \K manifolds is necessarily holomorphic or antiholomorphic, see Lemma~\ref{LemIsometries} below.
Notice also that weak relatives manifolds are genuinely relatives when the manifolds $X_1$ and $X_2$ in the definition are in fact locally biholomorphically isometric. 
This is not always the case.
For instance, it suffices to consider a $2$-dimensional \K manifold $(M,g)$ with discrete isometry group $G$ and a continuous family of parallel complex structures. 
If $J_1$ and $J_2$ are two complex structure which lie in different $G$-orbits, then the \K manifolds $(M,g,J_1)$ and $(M,g,J_2)$ are weak relatives but not relatives.
For example, this is the case for hyperK\"ahler manifolds. 
For dimension larger than $2$ it is not clear whether there exist weak relatives manifolds which are not relatives.

Our main theorem asserts that being weak relatives is not in fact a weaker condition than being relatives when one of the \K manifolds considered is projective, that is, is equipped with a metric induced by a holomorphic immersion in $\CP^N$.
	\begin{theor}\label{ThmMain}
	If a projective manifold $M_1$ and a \K manifold $M_2$ are weak relatives, then they are relatives.
\end{theor}

This result should be compared to \cite[Lemma~2.4 and Remark~2.5]{DiscalaLoi10Relatives} where the same conclusion is obtained when $M_2$ is a homogeneous bounded domain of nonpositive holomorphic bisectional curvature.
Here we have no assumptions on the metric of $M_2$ nor on its curvature.
Theorem~\ref{ThmMain} allows us to prove that certain pairs of \K manifolds are not weak relatives whenever we can show that they are not relatives and one of them is projective.
In particular, it was proven in \cite[Theorem~1.2]{LoiMossa23Rigidity} that a complex projective space and the product of a flat \K manifold with a homogeneous bounded domain are not relatives.
From this we can draw the following corollary.
\begin{cor}
A projective \K manifold $X$ and a product $\C^N/\Gamma\times\Omega$ of a flat \K manifold with a homogeneous bounded domain $\Omega$ are not weak relatives.
\end{cor}
For instance, a complex torus with a projective \K metric and a flat complex torus are not weak relatives. 
Notice that this cannot be deduced without Theorem~\ref{ThmMain}, although it was already known to Umehara \cite{Umehara87Relatives} that they are not relatives.

Notice that, if there exists a holomorphic isometric immersion $M_1\lra M_2$ of \K manifolds, then trivially $M_1$ and $M_2$ are relatives.
However, the relation of being relatives does appear to be of local nature. 
In fact, as observed after Definition~\ref{DefRelatives}, one can always restrict to an open set where the immersions are embeddings. 
Even more so, in many cases the obstruction to being relatives can be detected in any open set of a complex curve $U\subset\C\lra M_i$ around a point, see for instance \cite{LoiMossa23Rigidity}. 
In light of this, it is surprising that the only examples of related \K manifolds $M_1$, $M_2$ in the literature are due to the existence of a holomorphic isometry $M_1\lra M_2$, or viceversa. 
That is, the known examples rely on a solution of the original problem of Calabi and not strictly on its generalization. 
We find that studying related \K manifolds which cannot be immersed into one another could give a better understanding of this relation and how it differs from Calabi's original problem. 
For this reason we propose the following
\begin{defin}\label{DefStrictlyRelatives}
Two \K manifolds $M_1$ and $M_2$ are \textbf{strict relatives} (or are strictly related) if they are relatives and there exists no local holomorphic isometry of one into the other.
\end{defin}
Observe that, if a \K manifold is related to a product of \K manifolds, then it is related to its factors.
Therefore, it is rather easy to produce examples in which one of the manifolds considered is a product, cf. Example~\ref{ExNonIrriducible}. 
Considering this, we provide several examples of strictly related irreducible \K manifolds.
Namely, Example~\ref{ExIrriducibleNonCpt} and Example~\ref{ExIrriducibleNonFlat} exhibits two pairs of irreducible noncompact \K manifolds  with different curvature assumptions.
In Example~\ref{ExIrriducibleCpt} we produce an instance in which both \K manifolds are compact.
Finally, in Example~\ref{ExIrriducibleCptNonCpt} we consider a pair where both \K manifolds are irreducible but only one is compact.
These are collected in Section~\ref{SecExamples} below.
	
	\section{Weak relatives are relatives}\label{sectionmain}
The aim of this section is to prove Theorem~\ref{ThmMain}. 
In order to do so we recall a result on isometries of \K manifolds. Although this is a well known result, to the best of the author's knowledge, there is no complete proof in the literature. Therefore we include here a simple proof for the reader's convenience.
\begin{lem}\label{LemIsometries}
Let $\varphi:M_1\lra M_2$ be an isometry between irreducible \K manifolds $M_1$ and $M_2$. If $M_1$is not Ricci-flat, then $\varphi$ is either holomorphic or anti-holomorphic.
\end{lem}
\begin{proof}
Denote by $J_1$ (respectively $J_2$) the complex structure on $M_1$ (resp. $M_2$).
Let $A\vcentcolon= \varphi^*(J_2)$ be the endomorphism of $TM_1$ given by the pullback of the complex structure of $M_2$.
Since $\varphi$ is an isometry and $M_2$ is K\"ahler, we have $\nabla A=0$.
Thus, for all $p\in M_1$, $A_p:T_pM_1\lra T_pM_1$ commutes with the Lie algebra $\hol^o_p(M_1)$ of the restricted holonomy group $\Hol^o_p(M_1)$ of the metric $g_1$ at $p\in M_1$.
Clearly, the same holds true for $J_1$.
Now a result due to Lichnerowicz (see for instance \cite[Proposition~4.1]{DiScala02Lichnerowicz}) ensures that the Lie algebra $\hol^o_p(M_1)$ of $\Hol^o_p(M_1)$ coincides with its normalizer in the Lie algebra $\so(T_pM_1)$ since $M_1$ is irreducible and not Ricci-flat.
As a consequence, $A_p$ and $J_1$ belong to the Lie algebra $\hol^o_p(M_1)$.
Therefore, we can identify $A_p$ with an automorphism of $\C^n\cong T^{1,0}_pM_1$ such that $J_1$ is identified with multiplication by $i$.
By irreducibility of $M_1$, the Lie algebra $\hol^o_p(M_1)$ acts irreducibly on $T_pM_1$.
Therefore Schur's lemma implies that $A_p=\lambda\Id$ for some $\lambda\in\C$.
Thus we have $\lambda^2\Id=A_p^2=-\Id$.
We conclude $\lambda=\pm i$ so that $A_p=\pm J_1$ for all $p\in M_1$.
In other words, $A= \varphi^*(J_2)=\pm J_1$ as wanted.
\end{proof}

\begin{proof}[Proof of Theorem~\ref{ThmMain}]
Let $M_1$ be a projective manifold, that is, $M_1$ admits a holomorphic immersion $\iota:M_1\lra\CP^N$ and is equipped with the pullback of the Fubini-Study metric via such immersion.
Suppose now the \K $M_2$ is a weak relative of $M_1$.
Namely, there exists locally isometric \K manifolds $X_1$ and $X_2$ and holomorphic isometries $\varphi_i:X_i\lra M_i$ for $i=1,2$.
Denote by $\varphi$ a local isometry between $X_1$ and $X_2$. 
Since we are only interested in the local behaviour, without loss of generality, we can assume that $\varphi$ is a global isometry between $X_1$ and $X_2$.
We want to show that we can modify $\varphi$ so to obtain a holomorphic isometry $\widetilde{\varphi}:X_1\lra X_2$ which implies that  $M_1$ and $M_2$ are relatives. 

Let $X_1=F\times N_1\times\cdots\times N_k$ be the de Rham decomposition of $X_1$ where we collected in $F$ all the Ricci-flat factors.
By a result of Hulin \cite{Hulin96RicciFlat} (see \cite{ArezzoLiLoi23NonNegativeRicci} for a recent generalization) a Ricci-flat \K manifold does not admit a holomorphic isometry into a complex projective space. Therefore the factor $F$ cannot appear.
Now $\varphi$ restricts to each factor $N_j$ as an isometry of \K manifolds because $\varphi(N_j)$ must be an irreducible factor in the de Rham decomposition of $X_2$.
Denote this restriction by $\varphi_j\vcentcolon=\varphi_{\vert_{N_j}}$.
 By Lemma~\ref{LemIsometries},  $\varphi_j$ is either holomorphic or anti-holomorphic.
Now, let
$$
\widetilde\varphi_j=
\begin{cases}
\varphi_j& \mbox{if } \varphi_j \mbox{ is holomorphic;}\\
	\overline{\varphi}_j & \mbox{if } \varphi_j \mbox{ is anti-holomorphic}.
\end{cases}
$$
where $\overline{\varphi}_j$ denotes the conjugate of $\varphi_j$.
We conclude that $\widetilde{\varphi}\vcentcolon=\widetilde\varphi_1\times\cdots\times\widetilde\varphi_k$ is a holomorphic isometry between $X_1$ and $X_2$ and this concludes the proof.
\end{proof}

\section{Stricltly related \K manifolds}\label{SecExamples}
As mentioned in Section~\ref{sectionint}, there has been a lack of examples of strict relatives \K manifolds in the literature.
In this section we provide several new examples of \K manifolds $(M_1,g_1)$ and $(M_2,g_2)$ which are strict relatives, that is, which are related but such that one cannot be locally holomorphically isometrically immersed into the other.
We begin by the simplest case and later provide gradually more restrictive examples.
\begin{ex}\label{ExNonIrriducible}\rm
Consider $M_1=\C^n$ endowed with the flat metric and the \K product $M_2=\C\times\CP^m$ so that $g_2$ is the sum of the standard flat metric with the Fubini-Study metric.
Clearly, the complex line $\C$ (endowed with the flat metric) can be embedded in both $M_1$ and $M_2$.

Notice that $M_2$ does not admit a local holomorphic isometry in $M_1$ because of a classical result of Calabi \cite{Calabi53Isometric} on holomorphic isometries of complex space forms.
On the other hand, if we assume $ n>m+1 \geq 2$, then $M_1$ cannot be holomorphically isometrically immersed in $M_2$ for dimensional reasons. 
Moreover, $M_1$ cannot be holomorphically isometrically immersed in $M_2$ even for $n=m+1$ because such an immersion would be a local isometry which is obstructed by the holomorphic sectional curvature.
Hence, for $n\geq m+1\geq 2$, $M_1$ and $M_2$ are strict relatives.
\end{ex}

In the previous example, the common \K submanifold $\C$ is a \K factor of both $M_1$ and $M_2$ so it trivially embeds in them.  
It is natural to ask whether there exist strictly related irreducible manifolds. 
In the following example we provide one such instance where the manifolds are non-compact .

\begin{ex}\label{ExIrriducibleNonCpt}\rm
Let $M_1=\C^k\#\overline{\CP}^k$ be the blow-up of $\C^k$ at the origin endowed with a multiple $\eta g_S$ of the Burns-Simanca metric $g_{S}$ such that $0<\eta\in\R$.
Recall that this metric is scalar flat, see for instance \cite{CannasAghedu19Simanca,LeBrun88SimancaMetric}.
Moreover, consider a homogeneous \K manifold $M_2=\C^n\rtimes\CH^m$ being the complex (but not K\"ahler) product of $\C^n$ with the disc $\CH^m$ endowed with a large enough multiple $\lambda g$ of a homogeneous \K metric $g$, see for example \cite{Yang07TwistProduct}.
Now the complex line $\C$ (endowed with the flat metric) can be embedded in both $M_1$ (see \cite[Theorem~1.1]{LoiMossa22BlowUp}) and $M_2$ (in a fiber of the holomorphic fibration, i.e. $z\mapsto(z,z_2,\ldots,z_n,x)$).

It was proven in \cite{Loi15HomogeneousImmersions} that $M_2$, being a simply connected homogeneous \K manifold, can be holomorphically isometrically immersed into $\CP^\infty$ when endowed with a large enough real multiple $\lambda g$ of any homogeneous \K metric $g$. We can then assume that $\lambda$ is such.
Furthermore, an open set $U$ in $M_1$ admits a holomorphic isometric immersion into $\CP^\infty$ only if $\eta$ is a positive integer. This follows from the fact that any such immersion of $U$ into $\CP^\infty$ can be extended to a global holomorphic isometric immersion of $M_1$ by Calabi's extension theorem \cite{Calabi53Isometric}, being $M_1$ simply connected.  On the other hand, any projectively induced metric must be integral. 
This only happens if $\eta$ is an integer because $0\neq[\omega_S]\in H^2(M_1;\Z)$ where $\omega_S$ is the \K form associated to $g_S$.

We can now assume that $n+m>k$ to prevent the existence of an immersion $M_2\lra M_1$.
In order to rule out the existence of an immersion of $M_1$ into $M_2$ it is enough to choose a large enough $\lambda\in\R\setminus\Z$.
Notice that this argument also excludes the case $n+m=k$ because in that case an immersion is a local isometry.
\end{ex}

We give now an example in which the manifolds $M_1$ and $M_2$ are still irreducible and noncompact as in Example~\ref{ExIrriducibleNonCpt} but not scalar flat as in the case of the Burns-Simanca metric.
\begin{ex}\label{ExIrriducibleNonFlat}\rm
	Let $M_1=\CH^n$ with the hyperbolic metric and let $M_2$ be a bounded symmetric domain of dimension $m\leq m$ and rank$\geq2$ equipped with the Bergman metric.
	Both $M_1$ and $M_2$ admit $\CH^1$ as a totally geodesic submanifold.
	
	Clearly, if $n>m$, $M_1$ does not admit a local immersion in $M_2$ for dimensional reasons. As in the previous example, this holds also when $n=m$ because a local holomorphic isometry $M_2\lra M_1$ would force $M_2$ to have constant holomorphic sectional curvature. 
	On the other hand, if $M_2$ admitted a local holomorphic isometric immersion into $M_1$, then it would admit one into the infinite dimensional simply connected space form $\ell^2(\C)$ by composing with the immersion $\CH^n\lra\ell^2(\C)$. However, by \cite[Theorem~3.3]{DiScalaLoi07SymmetricSpaces} this is only possible when $M_2$ is (a \K product of copies of) $\CH^k$.
	Therefore, for $n\geq m$, $M_1$ and $M_2$ are strictly related.
\end{ex}

One would expect that compactness does not play a role in a problem of local nature as finding non strict relatives. Indeed, our next examples exhibits two compact \K manifold which are strict relatives.

\begin{ex}\label{ExIrriducibleCpt}\rm
	Let $M_1=\CP^n$ with the Fubini-Study metric and let $M_2=Q^m$ where $Q^m$ is the $m$-dimensional quadric in $\CP^{m+1}$ with the projective homogeneous metric.
	Both $M_1$ and $M_2$ admit $\CP^1$ as a totally geodesic submanifold. 
	
	If we assume $n<m$, $M_2$ does not admit a local immersion in $M_1$ for dimensional reasons.
	As in the previous example, also the case $n=m$ is ruled out because the existence of a local holomorphic isometry $M_2\lra M_1$ would imply that $M_2$ has constant holomorphic sectional curvature.
	Viceversa, by \cite[Theorem~1]{Suyama82ProjectiveIntoQuadrics} $M_1$ can be (locally) holomorphically isometrcally immersed into $M_2$ only if $m\geq 2n$. Thus, for $n\leq m<2n$, $M_1$ and $M_2$ are strict relatives.
\end{ex}

We conclude this section by showing an example of two strictly related irreducible \K manifolds such that only one is compact.

\begin{ex}\label{ExIrriducibleCptNonCpt}\rm
	Let $M_1=\CP^n$ with the Fubini-Study metric. 
	The definition of the manifold $M_2$ is slightly more involved. Namely, $M_2$ as a complex manifold is the projectivization of the tangent bundle of $\CH^2$. This can be seen as the noncompact homogeneous \K manifold $SU(2,1)/U(1)\times U(1)$ with the natural $SU(2,1)$-action on the tangent bundle of $\CH^2$.
	Let $z_1,z_2$ be coordinates on $\CH^2$ and $w_0,w_1$ be homogeneous coordinates on the fiber $\CP^1$.
	Then consider the homogeneous metric $g_2$ on $M_2$ given in the dense open set $\{w_0\neq 0\}$ by 
	$$g_2=-2i\de\deb\log(1-\vert z_1\vert^2-\vert z_2\vert^2)+i\de\deb\log(1+\vert w\vert^2-\vert z_1+ z_2w\vert^2)$$
	where $w=w_1/w_0$. Notice that $M_2$ is the complex product of $\CP^1$ and $\CH^2$ being a bundle over a contractible space but the \K metric is not the product metric. 
	
	One easily sees that $\CP^1$ embeds into $M_2$ as the fiber over $z_1=z_2=0$ because  the metric restricted to that fiber reads $g_{2\vert_{\CP^1}}=\frac{i}{2}\de\deb\log(1+\vert w\vert^2)^2= \omega_{FS}$.
	Thus both $M_1$ and $M_2$ admit $\CP^1$ as a \K submanifold, the embedding $\CP^1\lra M_1$ being the totally geodesic one.
	If we assume $n\geq 3$, $M_1$ does not admit a local holomorphic isometric immersion in $M_2$ for dimensional reasons.
	In particular, the case $n=3$ is ruled out because the existence of a local holomorphic isometry $M_1\lra M_2$ would imply that $M_2$ has constant holomorphic sectional curvature in a neighbourhood of a point.
	
	On the other hand, one can show \cite{Zuddas23PrivateComm} that the metric $g_2$ is induced by the full immersion $F:M_2\lra\CP^\infty$ given in the coordinates above by
	$$ F(z_1,z_2,w)=\left(1,w,\ldots, \sqrt{{{k+1}\choose{j+1}} (k-j+1)} z_1^j z_2^{k-j}\left(  z_1-\dfrac{j+1}{k-j+1} z_2w\right),\ldots\right)$$
	for $k\geq 0$ and $0\leq j\leq k$.
	This prevents the existence of a local holomorphic isometry $U\subset M_2\lra M_1$. 
	In fact, if such an immersion existed, we would get two local holomorphic isometries of $M_2$ into $\CP^\infty$. 
	By Calabi's rigidity these are related by a rigid transformation of $\CP^\infty$ which is impossible because one isometry is full and the other is contained in a proper totally geodesic submanifold.
	Therefore, for $n\geq 3$, $M_1$ and $M_2$ are strict relatives.
\end{ex}

	\bibliographystyle{amsplain}
	
	\bibliography{/Users/giovanniplacini/Library/CloudStorage/Dropbox/Projects/RTT_PROGETTI/Bibliography/biblio.bib}

\end{document}